\documentclass[12pt]{article}

\usepackage{amsfonts,graphicx}

\newtheorem{lemma}{Lemma}
\newtheorem{theorem}{Theorem}

\title{On convexity of the frequency response of\\
a stable polynomial}

\author{Didier Henrion$^{1,2}$}

\begin{document}

\maketitle

\footnotetext[1]{LAAS-CNRS, University of Toulouse, France}
\footnotetext[2]{Faculty of Electrical Engineering,
Czech Technical University in Prague, Czech Republic}

\begin{abstract}
In the complex plane, the frequency response of a univariate polynomial
is the set of values taken by the polynomial when evaluated
along the imaginary axis. This is an algebraic curve
partitioning the plane into several connected components.
In this note it is shown that the component including the origin
is exactly representable by a linear matrix
inequality if and only if the polynomial is stable,
in the sense that all its roots have negative real parts.
\end{abstract}

\begin{center}
{\bf\small Keywords}\\[1em]
Polynomial, stability, convexity, linear matrix inequality,
real algebraic geometry.\\
\end{center}

\section{Introduction}

Let $p(s) \in {\mathbb R}[s]$
be a polynomial of the complex indeterminate $s \in {\mathbb C}$.
We say that $p(s)$ is stable if all its roots lie in the
open left half-plane.
Define the frequency response
\[
{\mathcal P} = \{p(j\omega) \: :\: \omega \in {\mathbb R}\}
\]
as the set of values taken by the polynomial when evaluated
along the stability boundary, namely the imaginary axis.
The frequency response plays
a key role when deriving results of robust control theory
such as Kharitonov's theorem \cite{ackermann,barmish,bhatta}.

In \cite{hamann} it was observed that the frequency response
of a stable polynomial features interesting convexity properties,
see also \cite[Chapter 18]{barmish}.
More specifically, given a polynomial $p(s)$, 
an arc is defined as a subset of the frequency response
for a given range of the indeterminate, i.e.
$\{p(j\omega) \: :\: \omega \in \Omega\}$ with
$\Omega$ a subset of the real line $\mathbb R$.
A proper arc is an arc that does not pass through the
origin and such that the net change in the argument
of $p(j\omega)$ does not exceed $\pi$ as $\omega$
increases over $\Omega$. The Arc Convexity Theorem of
\cite{hamann}
states that all proper arcs of the value set of a
stable polynomial are convex. Alternative proofs
can be found in \cite{kogan} and \cite{gu}.

The frequency response $\partial{\mathcal P}$ is a curve
that partitions the complex plane into several connected
components. We denote by $\mathcal P$ the connected
component including the origin. It is called the
inner frequency response set in \cite{hamann}.
In the case of a stable polynomial, the boundary
of $\mathcal P$ therefore consists of a finite union of
proper arcs. Theorem 4.1 in \cite{hamann} uses the
Arc Convexity Theorem to establish convexity of
$\mathcal P$.

In this note we provide a more accurate description
of the geometry of $\mathcal P$ and an alternative proof
of its convexity. We derive an explicit representation
of this set as a two-dimensional linear matrix inequality
(LMI).
Instrumental to this derivation are standard results
from real algebraic geometry \cite{abhyankar,roy,mourrain,sturmfels}
and a recent characterization of two-dimensional convex
polynomial level sets obtained in \cite{hv}.

\section{Algebraic description of the frequency response}

The frequency response of polynomial $p(s) = p_0+p_1s+\cdots+p_ns^n$
can be expressed as the parametric curve
\begin{equation}\label{param1}
\partial{\mathcal P} = \{p(j\omega) = x(\omega)+jy(\omega) \: :\:
\omega \in \mathbb R\}
\end{equation}
where
\begin{equation}\label{param2}
\begin{array}{rcl}
x(\omega) & = & \displaystyle{\frac{q_x(\omega)}{q_z(\omega)}}\\
y(\omega) & = & \displaystyle{\frac{q_y(\omega)}{q_z(\omega)}}\\
\end{array}
\end{equation}
and
\[
\begin{array}{rcl}
q_x(\omega) & = & p_0 - p_2 \omega^2 + p_4 \omega^4 + \cdots \\
q_y(\omega) & = & p_1\omega - p_3 \omega^3 + p_5 \omega^5 + \cdots \\
q_z(\omega) & = & 1
\end{array}
\]
are polynomials of the real indeterminate $\omega \in \mathbb R$.

Curve $\partial{\mathcal P}$ is rationally (here polynomially)
parametrized, so it is an algebraic plane curve of genus
zero \cite{abhyankar}. In control theory terminology,
$\partial{\mathcal P}$ is sometimes called the Mikhailov plot of
polynomial $p(s)$ or the Nyquist plot of
the rational (here polynomial) transfer function $p(s)$, see
\cite{barmish} or \cite{bhatta}.

Equations (\ref{param1}-\ref{param2}) provide a parametric description
of curve $\partial{\mathcal P}$. With the help of elimination
theory and resultants, we can derive an implicit description
\begin{equation}\label{impli}
\partial{\mathcal P} = \{x+jy \: :\: f(x,y) = 0\}
\end{equation}
where $f(x,y)$ is an irreducible bivariate polynomial, see
\cite[Section 3.3]{cox}.

\begin{lemma}\label{resultant}
Given two univariate polynomials $g(\omega) = g_0+g_1\omega+\cdots+g_n\omega^n \in
{\mathbb R}[\omega]$ and $h(\omega) = h_0+h_1\omega+\cdots+h_n\omega^n \in
{\mathbb R}[\omega]$, there exists a unique (up to sign) irreducible
polynomial $b(g,h) \in {\mathbb R}[g_0,g_1,\ldots,g_n,h_0,h_1,\ldots,h_n]$
called the resultant which vanishes whenever $g(\omega)$ and $h(\omega)$ have
a common zero.
\end{lemma}

To address the implicitization problem,
we make use of a particular resultant, the B\'ezoutian, see e.g.
\cite[Theorem 4.1]{sturmfels} or \cite[Section 5.1.2]{mourrain}.
Given two univariate polynomials $g(\omega), h(\omega)$ of the same degree $n$
as in Lemma \ref{resultant}, build the following bivariate polynomial
\[
\frac{g(\omega)h(v)-g(v)h(\omega)}{\omega-v} = \sum_{k=0}^{n-1} \sum_{l=0}^{n-1}
b_{kl} \omega^k v^l
\]
called the B\'ezoutian of $g$ and $h$,
and the corresponding symmetric matrix $B(g,h)$
of size $n\times n$ with entries $b_{kl}$ bilinear in coefficients
of $g$ and $h$. As shown e.g. in \cite[Theorem 4.3]{sturmfels} or
\cite[Section 5.1.2]{mourrain},
the determinant of the B\'ezoutian matrix is the resultant:

\begin{lemma}\label{bezout}
$\det B(g,h) = b(g,h)$.
\end{lemma}

Now we can use the B\'ezoutian to derive the
implicit equation (\ref{impli}) of curve $\partial{\mathcal P}$
from the explicit equations (\ref{param1}-\ref{param2}).

\begin{lemma}\label{pencil}
Given polynomials $q_x, q_y, q_z$ of 
equations (\ref{param1}-\ref{param2}), the polynomial
of equation (\ref{impli}) is given by
$f(x,y) = \det F(x,y)$ where
\begin{equation}\label{F}
F(x,y) = \pm(B(q_x,q_y)-xB(q_y,q_z)-yB(q_x,q_z))
\end{equation}
is a symmetric pencil, i.e. a polynomial matrix linear
in $x,y$, of size $n$.
\end{lemma}

{\bf Proof:} Rewrite the system of equations (\ref{param2}) as
\[
\begin{array}{rcccl}
g(\omega) & = & q_x(\omega) - x q_z(\omega) & = & 0 \\
h(\omega) & = & q_y(\omega) - y q_z(\omega) & = & 0 \\
\end{array}
\]
and use the B\'ezoutian resultant of Lemma \ref{bezout}
to eliminate indeterminate $\omega$ and obtain conditions
for a point $(x,y)$ to belong to the curve. The B\'ezoutian
matrix is $B(g,h) = B(q_x-xq_z,q_y-yq_z) =
B(q_x,q_z)-xB(q_y,q_z)-yB(q_x,q_z)$. Linearity in $x,y$
follows from bilinearity of the B\'ezoutian.
Finally, note that the sign of $F(x,y)=\pm B(g,h)$
affects the sign of $f$, but not the implicit
description $f(x,y)=0$.$\Box$

Lemma \ref{pencil} provides the implicit equation of
curve $\partial{\mathcal P}$ in symmetric linear determinantal form.

\section{Convexity properties of the inner frequency response set}

Curve $\partial{\mathcal P}$ partitions the complex plane into several
connected regions. We are interested in the connected region containing
the origin, denoted by $\mathcal P$.
In order to study the geometry of this region,
we need the following result.

\begin{lemma}\label{definite}
The sign of pencil $F(x,y)$ in (\ref{F}) can be chosen such that
$F(0,0)=B(q_y,q_z)$ is positive definite
if and only if $p$ is a stable
polynomial.
\end{lemma}

{\bf Proof:} The signature of the B\'ezoutian matrix $B(q_x,q_y)$,
(the number of positive eigenvalues minus the number
of negative eigenvalues) is equal to the Cauchy index of
the rational function $q_x(\omega)/q_y(\omega)$ (the number
of jumps from $-\infty$ to $+\infty$ minus the number of
jumps from $+\infty$ to $-\infty$), see \cite[Section 9.1.2]{roy}. 
The Cauchy index is maximum (resp. minimum) when
$B(q_x,q_z)$ is positive (resp. negative) definite. This
occurs if and only if polynomials $q_x$ and $q_y$ satisfy the root
interlacing condition, i.e. they must have only real roots and
between two roots of $q_x(\omega)$ there is only one root of
$q_y(\omega)$ and vice-versa. Since $q_x(\omega) = \mathrm{Re}\:p(j\omega)$
and $q_y(\omega) = \mathrm{Im}\:p(j\omega)$, this is equivalent
to stability of $p$ in virtue of the Hermite-Biehler
theorem, see \cite[Section 8.1]{ackermann} or
\cite[Section 1.3]{bhatta}.$\Box$

The main result of this note can now be stated.

\begin{theorem}\label{convex}
The connected component including the origin
and delimited by the frequency response of polynomial $p$
can be described by a linear matrix inequality (LMI)
\[
{\mathcal P} = \{x+jy \: :\: F(x,y) \succeq 0\}
\]
if and only if $p$ is stable.
In the above description, $F(x,y)$ is given by (\ref{F})
and $\succeq 0$ means positive semidefinite.
\end{theorem}

{\bf Proof:}
First we prove that $p$ stable implies LMI 
representability of $\mathcal P$.
This set is the closure of
the connected component of the polynomial level set
$\{x,y \: :\: f(x,y) > 0\}$
that contains the origin, an algebraic interior
in the terminology of \cite{hv}.
Polynomial $f(x,y)$ is called the defining polynomial.
By continuity, the boundary of $\mathcal P$ consists
of those points $x+jy$ for which $F(x,y)$ drops rank
while staying positive semidefinite. Note that in general
this boundary is only a subset of the curve
$\partial{\mathcal P}$.

To prove the converse, namely that
LMI representability of $\mathcal P$
implies stability of $p$, we use a result
of \cite{hv} stating
that a two-dimensional algebraic interior containing
the origin has an LMI representation if and only if it is
rigidly convex. Geometrically, this means that a generic
line passing through the origin must intersect the
algebraic curve $f(x,y) = 0$ a number of times equal
to the degree of $f(x,y)$. Rigid convexity implies that
$x(\omega)$ and $y(\omega)$, the respective real and
imaginary parts of polynomial $p(j\omega)$, satisfy
the root interlacing property, and this implies
stability of $p$ by the Hermite-Biehler
Theorem used already in the proof of Lemma \ref{definite}.
$\Box$

\section{Examples}

\subsection{Stable third degree polynomial}

\begin{figure}
\begin{center}
\includegraphics[scale=0.8]{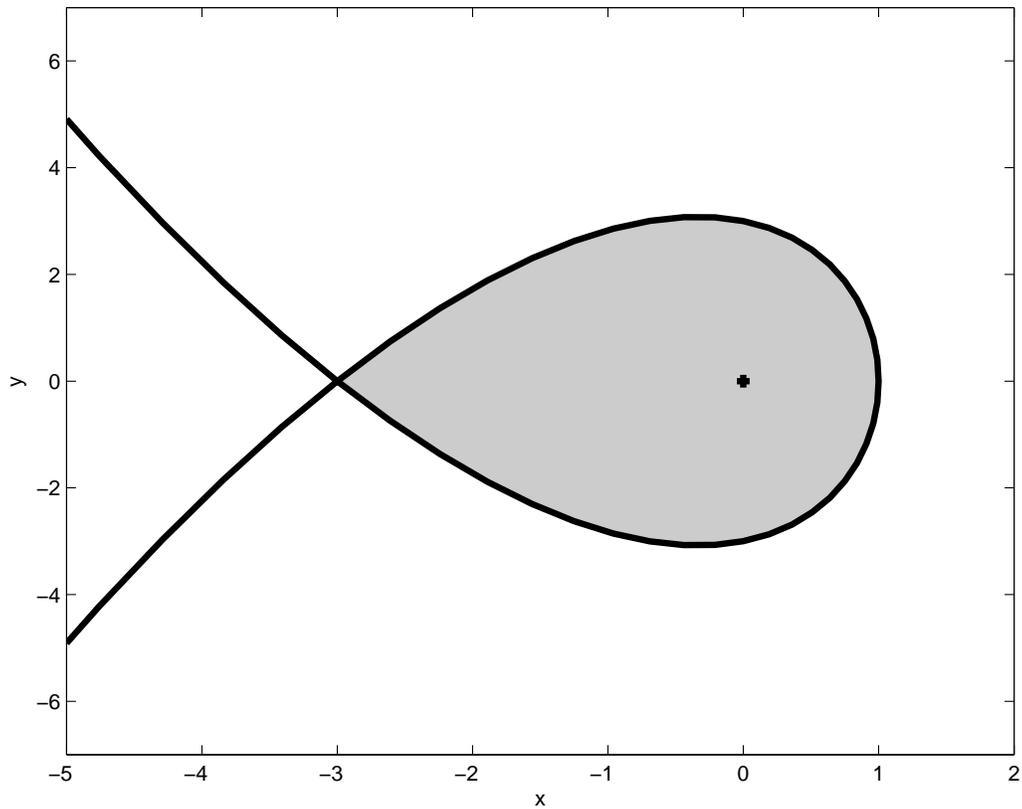}
\caption{Frequency response of a stable third degree polynomial. The shaded
region is the convex component containing the origin.\label{stable3}}
\end{center}
\end{figure}

Let $p(s)=s^3+s^2+4s+1$. Then $q_x(\omega) = 1-\omega^2$
and $q_y(\omega) = 4\omega-\omega^3$ in parametrization
(\ref{param2}). With the help of the Control System Toolbox
for Matlab, a visual representation of curve $\partial{\mathcal P}$
can be obtained as follows:
\begin{verbatim}
>> p = [1 1 4 1]; % polynomial in Matlab format
>> nyquist(tf(p,1)) % frequency response
>> axis([-5 2 -7 7]) % zoom around the origin
\end{verbatim}
see Figure \ref{stable3}.
The B\'ezoutian matrices of Lemma \ref{pencil} can be computed 
with the following Maple 10 instructions:
\begin{verbatim}
> with(LinearAlgebra):
> qx:=1-w^2:qy:=4*w-w^3:qz:=1:
> Bxy:=BezoutMatrix(qx,qy,w,method=symmetric);
                                   [-1     0     1]
                                   [              ]
                            Bxy := [ 0    -3     0]
                                   [              ]
                                   [ 1     0    -4]

> Byz:=BezoutMatrix(qy,qz,w,method=symmetric);
                                   [ 0     0     1]
                                   [              ]
                            Byz := [ 0     1     0]
                                   [              ]
                                   [ 1     0    -4]
> Bxz:=subs(e=0,BezoutMatrix(qx,qz+e*w^3,w,method=symmetric));
                                   [ 0     0     0]
                                   [              ]
                            Bxz := [ 0     0    -1]
                                   [              ]
                                   [ 0    -1     0]
\end{verbatim}
Note the use of the {\tt subs} instruction to ensure that
the last B\'ezoutian matrix has appropriate dimension 3.
Matrix $B(q_x,q_z)$ is negative definite,
so a sign change is required to build
\[
F(x,y) = \left[\begin{array}{ccc}
1 & 0 & -1+x \\ 0 & 3+x & -y \\ -1+x & -y & 4-4x
\end{array}\right]
\]
and we obtain the following determinantal polynomial:
\begin{verbatim}
> F:=-(Bxy-x*Byz-y*Bxz);
> f:=Determinant(F);  
                                           2    2    3
                         f := 9 - 3 x - 5 x  - y  - x
\end{verbatim}
describing algebraic plane curve $\partial{\mathcal P}$ implicitly.
The curve can be studied with the {\tt algcurves} package
of Maple:
\begin{verbatim}
> with(algcurves):
> genus(f,x,y);
                         0
> plot_real_curve(f,x,y);
\end{verbatim}

\subsection{Stable eighth degree polynomial}

\begin{figure}
\begin{center}
\includegraphics[scale=0.8]{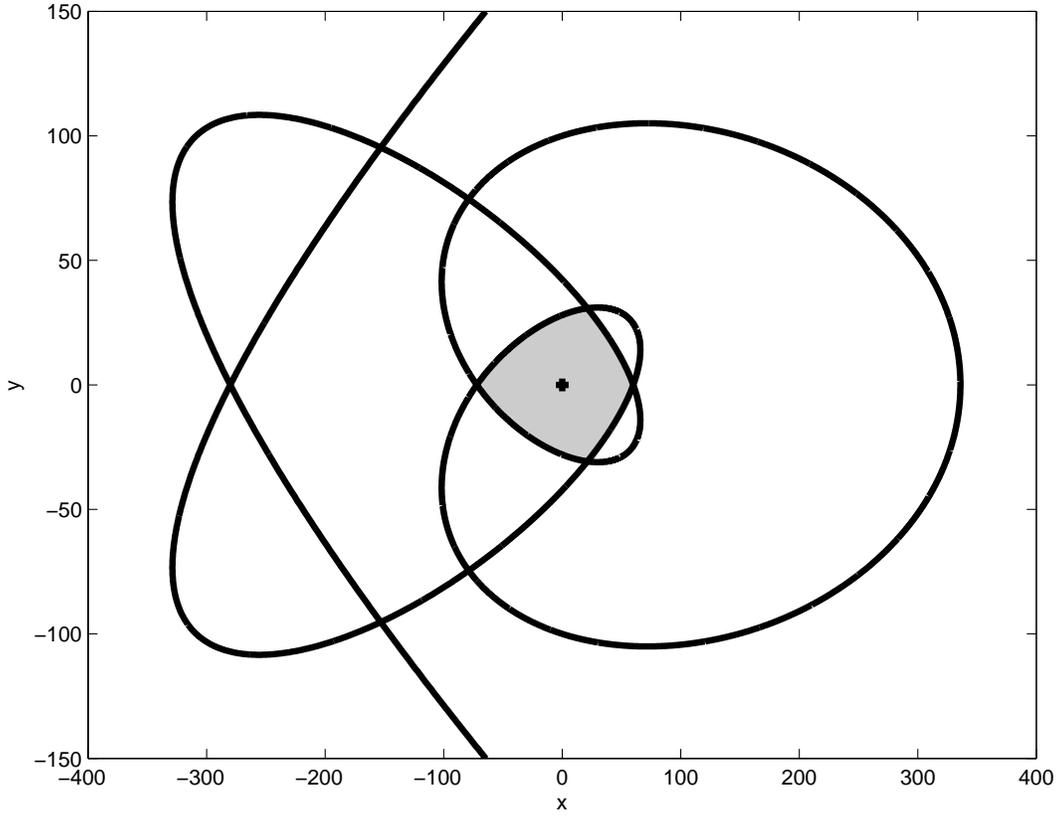}
\caption{Frequency response of a stable eighth degree polynomial.
The shaded region is rigidly convex.\label{stable8}}
\end{center}
\end{figure}

A more complicated example is the eighth degree stable polynomial
$p(s) = 336+198s+496s^2+117s^3+183s^4+20s^5+24s^6+s^7+s^8$
whose frequency response is represented on Figure \ref{stable8}.
The rigidly convex region around the origin has the LMI description

{\footnotesize\[
\left[ \begin {array}{c@{\;}c@{\;}c@{\;}c@{\;}c@{\;}c@{\;}c@{\;}c} 1&0&-20&0&117&0&-198&y\\
0&4&0&-66&0&298&y&-336+x\\
-20&0&414&0&-2510&y&4416+x&-24\,y\\
0&-66&0&1150&y&-5504+x&-24\,y&6720-20\,x\\
117&0&-2510&y&15907+x&-24\,y&-29514-20\,x&183\,y\\
0&298&y&-5504+x&-24\,y&28518-20\,x&183\,y&-39312+117\,x\\
-198&y&4416+x&-24\,y&-29514-20\,x&183\,y&58896+117\,x&-496\,y\\
y&-336+x&-24\,y&6720-20\,x&183\,y&-39312+117\,x&-496\,y&66528-198\,x
\end {array}\right] \succeq 0
\]}
\subsection{Unstable fourth degree polynomial}

\begin{figure}
\begin{center}
\includegraphics[scale=0.8]{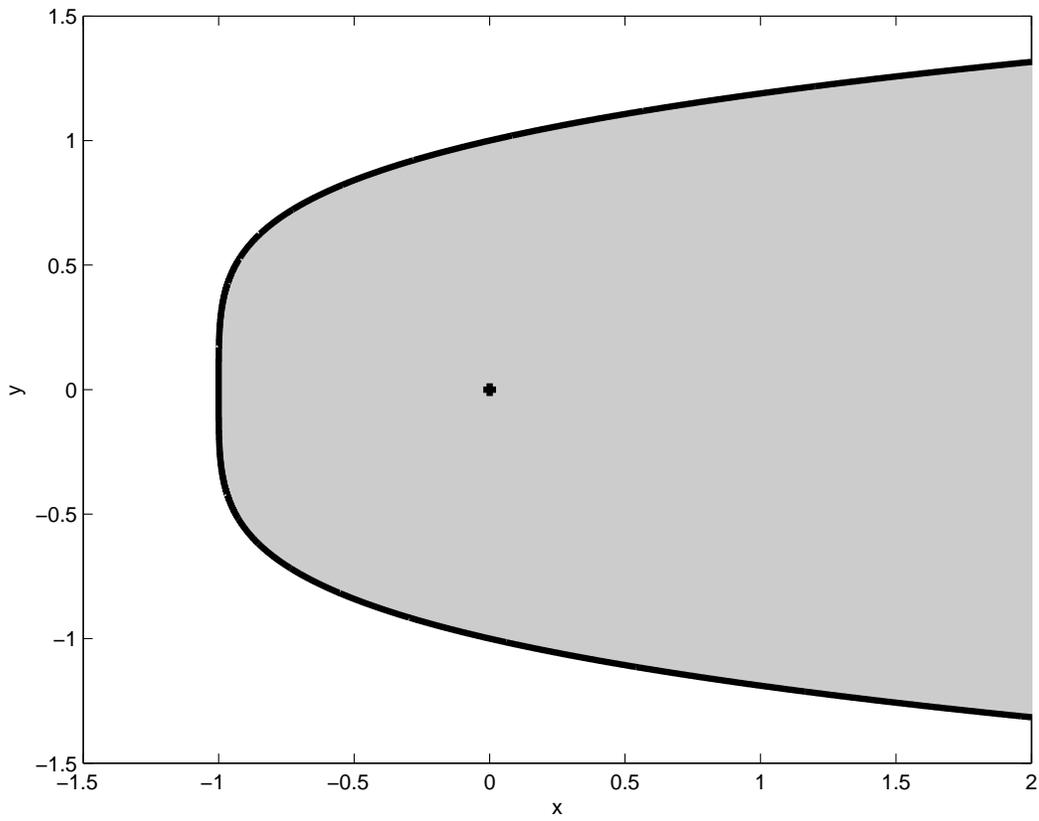}
\caption{Frequency response of an unstable fourth degree polynomial.
The shaded region is convex, but not rigidly convex.
\label{unstable4}}
\end{center}
\end{figure}

This example is taken from \cite{kogan}.
Let $p(s)=s^4-s-1$. The implicit equation of $\partial{\mathcal P}$
is
\[
f(x,y) = \det\:\left[\begin{array}{cccl}
0 & 0 & 1 & y\\
0 & 1 & y & 0\\
1 & y & 0 & 0\\
y & 0 & 0 & 1+x
\end{array}\right] =
-1-x+y^4 = 0.
\]
The component $\mathcal P$ including the origin
is convex, see Figure \ref{unstable4}, but it is not rigidly
convex since a generic line passing through the origin
cuts the quartic $\partial{\mathcal P}$ only twice. Hence this region
does not admit an LMI representation, and by Theorem \ref{convex},
polynomial $p(s)$ is unstable.

Note however that $\mathcal P$ can
be represented as the projection of an LMI set:
\[
{\mathcal P} = \{x+jy:\: \exists z:\:\left[\begin{array}{cc}
1+x & z\\
z & 1
\end{array}\right] \succeq 0, \:
\left[\begin{array}{cc}
z & y\\
y & 1
\end{array}\right] \succeq 0\}
\]
by introducing a lifting variable $z$, 
but such constructions are out of the scope of this note.

\section{Conclusion}

Convexity of the connected component
containing the origin and delimited by the frequency response
of a stable polynomial was already established in
\cite{hamann}. In this note we give an alternative
proof of this result based on B\'ezoutians
and we give a more accurate characterization of the geometry
of this region. Namely, the region is rigidly convex
in the sense of \cite{hv}, a property which is stronger than
convexity, and which is equivalent to the existence
of an LMI representation of the set.

In the terminology of convex analysis, the polynomial $f(x,y)$
defining implicitly the frequency response in (\ref{impli})
is hyperbolic with respect to the origin. Equivalently, $f(x,y)$
can be expressed as the determinant of a symmetric pencil which
is positive definite at the origin. See \cite{renegar}
for a tutorial on hyperbolic polynomials and \cite{lewis}
for connections with the results of \cite{hv}. This note
therefore unveils a link between polynomial
hyperbolicity and stability.

\section*{Acknowledgments}

This work benefited from discussions with Bernard Mourrain.

\end{document}